\documentclass[a4paper,12pt,final]{amsart}
\usepackage{times,a4wide,mathrsfs,amssymb,amsmath,amsthm,enumerate,xypic,tikzsymbols,dsfont}

\newcommand{\C}{\mathbb{C}}

\newcommand{\QQ}{\mathbb{Q}}
\newcommand{\NN}{\mathbb{N}}
\newcommand{\PP}{\mathbb{P}}

\newcommand{\OO}{\mathcal O}

\newcommand{\XX}{\mathcal X}
\newcommand{\YY}{\mathcal Y}

\newcommand{\MM}{\mathcal M}

\newcommand{\BB}{\mathfrak B}

\newcommand{\pic}{\hbox{Pic}}

\newcommand{\gr}{\hbox{Gr}}
\newcommand{\wt}{\widetilde}
\newcommand{\rom}{\romannumeral}

\newcommand{\one}{\mathds{1}}

\DeclareMathOperator{\aut}{Aut}

\DeclareMathOperator{\ide}{id}

\DeclareMathOperator{\ima}{Im}

\DeclareMathOperator{\Gr}{Gr}

\newtheorem{theorem}{Theorem}[section]

\newtheorem{lemma}[theorem]{Lemma}

\newtheorem{corollary}[theorem]{Corollary}
\newtheorem{proposition}[theorem]{Proposition}
\newtheorem{conjecture}[theorem]{Conjecture}
\newtheorem{remark}[theorem]{Remark}
\newtheorem{definition}[theorem]{Definition}
\newtheorem{convention}{Conventions}

\newtheorem{notation}[theorem]{Notation}

\newtheorem{nonumbering}{Theorem}

\newtheorem{nonumberingt}{Acknowledgements}

\begin{document}

\author[Robert Laterveer]
{Robert Laterveer}

\address{Institut de Recherche Math\'ematique Avanc\'ee,
CNRS -- Universit\'e 
de Strasbourg,\
7 Rue Ren\'e Des\-car\-tes, 67084 Strasbourg CEDEX,
FRANCE.}
\email{robert.laterveer@math.unistra.fr}

\title{The Beauville--Voisin conjecture for double EPW sextics}

\begin{abstract} We prove that the Beauville--Voisin conjecture is true for any double EPW sextic, i.e. the subalgebra of the Chow ring generated by divisors and Chern classes of the tangent bundle injects into cohomology.
 \end{abstract}

\thanks{Supported by ANR grant ANR-20-CE40-0023.}

%\keywords{Algebraic cycles, Chow group, motive, Beauville's ``splitting property'' conjecture, multiplicative Chow--K\"unneth decomposition, Fano threefold, homological projective duality}
%\subjclass[2010]{Primary 14C15, 14C25, 14C30.}

\maketitle

\section{Introduction}

Given a smooth projective variety $X$ over $\C$, let $A^i(X):=CH^i(X)_{\QQ}$ denote the Chow groups of $X$ (i.e. the groups of codimension $i$ algebraic cycles on $X$ with $\QQ$-coefficients, modulo rational equivalence). The intersection product defines a ring structure on $A^\ast(X)=\bigoplus_i A^i(X)$, the {\em Chow ring\/} of $X$ \cite{F}. 

For hyperk\"ahler varieties, it is famously expected that the Chow ring has particular behaviour:

\begin{conjecture}(Beauville--Voisin \cite{Beau3}, \cite{V17})\label{conj} Let $X$ be a hyperk\"ahler variety. The $\QQ$-subalgebra
  \[  \bigl\langle A^1(X), c_j(X)\bigr\rangle\ \ \subset\ A^\ast(X) \]
  (generated by divisors and Chern classes of the tangent bundle of $X$)
  injects into cohomology, under the cycle class map.
\end{conjecture}

This conjecture is true for K3 surfaces \cite{BV}. In dimension larger than 2, there is only one locally complete family for which the conjecture is currently known: the Fano varieties of lines on cubic fourfolds, for which Voisin has settled the conjecture by ingenious geometric arguments \cite{V17}.

In this paper, we establish Conjecture \ref{conj} for another locally complete family of hyperk\"ahler fourfolds: the double EPW sextics constructed and studied by O'Grady \cite{OG}, \cite{OG2}, \cite{OG3}, \cite{OG4}, \cite{OG5}. In fact, we prove something slightly stronger:

\begin{nonumbering}[=Theorem \ref{main}] Let $X$ be a smooth double EPW sextic. The $\QQ$-subalgebra
  \[  \bigl\langle A^1(X), c_j(X)\bigr\rangle\ \ \subset\ A^\ast(X) \]
  injects into cohomology, under the cycle class map.
  
  In addition, let $A^2(X)^+\subset A^2(X)$ denote the subgroup of cycles invariant under the covering involution $\iota$. The cycle class map induces injections
  \[  \bigl\langle A^1(X), c_j(X), A^2(X)^+\bigr\rangle\ \cap A^i(X)\ \hookrightarrow\ H^{2i}(X,\QQ) \ \ \ \hbox{for}\ i\ge 3\ .\]
  \end{nonumbering}
  
 Theorem \ref{main} is a ``standing on the shoulders of giants'' type of result; the argument relies in an essential way on earlier work by various authors. We use the results of Ferretti \cite{Fe} (who proved Conjecture \ref{conj} for double EPW sextics with Picard number 1) and Laterveer--Vial \cite{LV} (who proved the $i=4$ case of Theorem \ref{main}). The argument for one-cycles (i.e. the $i=3$ case of Theorem \ref{main}) closely follows the path laid out by Voisin for Fano varieties of lines on cubic fourfolds \cite{V17}. This is possible thanks to the work of Iliev--Manivel \cite{IM} (who provide a geometric construction of the general double EPW sextic in terms of conics on a Gushel--Mukai fourfold) and the work of Zhang \cite{Zh} (who established a quadratic relation for the variety of conics on a Gushel--Mukai fourfold). We also rely on work of Perry--Pertusi--Zhao \cite{PPZ}, who give an interpretation of the very general double EPW sextic as a moduli space; this allows to prove the ``Franchettina property'' (which is explained in subsection \ref{ss:gen} below) for the universal family of double EPW sextics.
 
 Following Voisin's pioneering work on Fano varieties of lines on cubic fourfolds \cite{V17}, we prove the $i=3$ case of Theorem \ref{main} by establishing an equality of correspondences using the quadratic relation of Zhang \cite{Zh}.
 The use of the Franchettina property in the argument implies that we obtain an equality of correspondences modulo {\em algebraic equivalence}. The Voevodsky--Voisin nilpotence theorem then allows to upgrade modulo rational equivalence. (Cf. Remarks \ref{pity} and \ref{rem} below for further comment on this detour via algebraic equivalence.)
   
 \vskip0.6cm

\begin{convention} In this paper, the word {\sl variety\/} will refer to a reduced irreducible scheme of finite type over $\C$. A {\sl subvariety\/} is a (possibly reducible) reduced subscheme which is equidimensional. 

{\bf All Chow groups will be with rational coefficients}: we will denote by $A_j(Y)$ the Chow group of $j$-dimensional cycles on $Y$ with $\QQ$-coefficients; for $Y$ smooth of dimension $n$ the notations $A_j(Y)$ and $A^{n-j}(Y)$ are used interchangeably. 
The notation $A^j_{hom}(Y)$ will be used to indicate the subgroup of homologically trivial cycles.
For a morphism $f\colon X\to Y$, we will write $\Gamma_f\in A_\ast(X\times Y)$ for the graph of $f$.

We will write $B^\ast(X)$ to indicate the groups of algebraic cycles with $\QQ$-coefficients modulo algebraic equivalence.

The contravariant category of Chow motives (i.e., pure motives with respect to rational equivalence as in \cite{Sc}, \cite{MNP}) will be denoted 
$\MM_{\rm rat}$.
%We will write $H^j(X)$ to indicate singular cohomology $H^j(X,\QQ)$.
\end{convention}

 \section{Preliminaries}

 \subsection{Double EPW sextics: the original construction}
 
 As the name suggests, double EPW sextics are double covers of certain sextic hypersurfaces:
	
	\begin{definition}[Eisenbud--Popescu--Walter \cite{EPW}] 
		Let $A\subset \wedge^3 \C^6$ be a subspace which
		is Lagrangian with respect to the symplectic form on $\wedge^3 \C^6$ given by
		the wedge product. The {\em EPW sextic associated to $A$\/} is
		\[ Z_A:= \Bigl\{  [v]\in \PP(\C^6)\ \vert\ \dim \bigl( A\cap ( v\wedge
		\wedge^2 \C^6)\bigr) \ge 1\Bigr\}\ \ \subset \PP(\C^6)\  .\]
		An {\em EPW sextic\/} is a hypersurface $Z_A$ for some $A\subset \wedge^3 \C^6$
		Lagrangian.
	\end{definition}

	\begin{theorem}[O'Grady]
		Let $Z$ be an EPW sextic such that the singular locus $S:=\hbox{Sing}(Z)$ is a
		smooth irreducible surface. Let 
		  \[ f\colon\ X\ \to\  Z\] 
		  be the double cover branched over
		$S$. Then $X$
		is a smooth hyperk\"ahler fourfold of K3$^{[2]}$-type (a so-called {\em double
			EPW sextic}), 
		and the class 
		  \[h:=f^*c_1(O_Z(1))\ \ \in \ A^1(X)\]
		  defines a polarization of square
		$2$ for the Beauville--Bogomolov form.
		Double EPW sextics form a $20$-dimensional locally complete family.
			\end{theorem}
	
	\begin{proof} 
		This is \cite[Theorem 1.1(2)]{OG2}. We remark that the hypothesis on
		$\hbox{Sing}(Z)$ is satisfied by a generic EPW sextic (more precisely, it
		suffices that the Lagrangian subspace $A$ be in $\hbox{LG}(\wedge^3 V)^0$, which
		is a certain open dense subset of $\hbox{LG}(\wedge^3 V)$ defined in
		\cite[Section 2]{OG2}). Letting $A$ vary in $\hbox{LG}(\wedge^3 V)^0$, one
		obtains a locally complete family with $20$ moduli (as observed in
		\cite[Introduction]{OG2}).		
	\end{proof}

 \subsection{The Chow ring of double EPW sextics}
 
 In this subsection, we recall what is currently known about the Beauville--Voisin conjecture for double EPW sextics. First, Ferretti has proven the conjecture for double EPW sextics with Picard number 1; more precisely:
 
 \begin{theorem}\label{fer}(Ferretti \cite{Fe0}, \cite{Fe}) Let $X$ be a double EPW sextic. Let $h\in A^1(X)$ be the polarization induced from the double cover, and let $c:=c_2(T_X)\in A^2(X)$.
 The $\QQ$-subalgebra
   \[ \langle h, c\rangle\ \ \subset\ A^\ast(X)\]
   injects into cohomology, under the cycle class map.
   \end{theorem}
   
   \begin{proof} This is \cite[Theorem 1.1]{Fe}. More details and background can be found in 
   \cite{Fe0}.
    \end{proof}
   
   With Charles Vial, we (partially reproved and) improved upon Ferretti's result:
   
   \begin{theorem}\label{lv}(Laterveer--Vial \cite{LV}) Let $X$ be a double EPW sextic. Let $A^2(X)^+\subset A^2(X)$ denote the subgroup of cycles invariant under the covering involution
   $\iota$. The cycle class map induces an isomorphism
  \[  \bigl\langle A^1(X), c_j(X), A^2(X)^+\bigr\rangle\ \cap A^4(X)\ \cong\ H^{8}(X,\QQ)\cong\QQ\ .\]
  \end{theorem}
  
  \begin{proof} This is \cite[Theorem 1]{LV}.
   \end{proof}
  
  \subsection{Double EPW sextics: the Iliev--Manivel construction}
  
  As recalled in this subsection, the general double EPW sextic can alternatively be constructed in terms of conics on a Gushel--Mukai fourfold.
  
  \begin{definition} An {\em ordinary Gushel--Mukai fourfold\/} is a smooth dimensionally transverse complete intersection
    \[ Y:= \Gr(2,5)\cap H\cap Q\ \ \ \subset\ \PP^9\ ,\]
    where $\Gr(2,5)$ is the Grassmannian of 2-dimensional subspaces of a fixed 5-dimensional complex vector space and $H, Q$ are a hyperplane resp. a quadric with respect to the Pl\"ucker embedding.
    \end{definition}
    
    Gushel--Mukai varieties have been studied in depth by Debarre and Kuznetzov \cite{D}, \cite{DK}, \cite{DK1}, \cite{DK3}, \cite{DK2}.

    \begin{notation} Let $F(G)$ denote the Hilbert scheme parametrizing conics contained in $\Gr(2,5)$. Given an ordinary Gushel--Mukai fourfold $Y$, let $F=F(Y)$ be the Hilbert scheme of conics contained in $Y$. We write $P\subset F\times Y$ for the universal conic, with projections $p\colon P \to F$ and $q\colon P\to Y$.
    
Conics in $\Gr(2,5)$ are of type $\tau$ or of type $\rho$ or of type $\sigma$ (cf. \cite[Section 3.1]{IM}); we write
    \[ F(Y)= F_\tau(Y)\cup F_\rho(Y)\cup F_\sigma(Y) \]
    for the partition of $F(Y)$ according to the type of conic.
      \end{notation}
      
      \begin{theorem}(Iliev--Manivel \cite{IM})\label{im} Let $Y$ be a general (ordinary) Gushel--Mukai fourfold. The variety $F=F(Y)$ is a smooth projective variety of dimension 5, and the subvarieties $F_\rho(Y), F_\sigma(Y)\subset F$ are of codimension 1 resp. 2. There exists a morphism
        \[ \pi\colon F\ \to\ X\ ,\]
        where $X$ is a double EPW sextic. The morphism $\pi$ restricted to $F_\tau(Y)$ is a $\PP^1$-fibration, and $\pi$ contracts $F_\rho(Y)$ and $F_\sigma(Y)$ to points $x_\rho$ resp. $x_\sigma$ on $X$.
        
        This construction defines a dominant rational map 
        \[\nu\colon\ \  \MM_{GM4}\ \dashrightarrow\ \MM_{dEPW} \]
        from the moduli space of Gushel--Mukai fourfolds to the moduli space of double EPW sextics.
        \end{theorem}
        
        \begin{proof} The smoothness of $F$ is \cite[Theorem 3.2]{IM}; the relation between $F$ and a double EPW sextic is \cite[Proposition 4.18]{IM}, and the fact that a general double EPW sextic is attained in this way is \cite[Corollary 4.17]{IM}.
         \end{proof}

\subsection{Double EPW sextics: the modular construction}    

More recently, it has been shown that the very general double EPW sextic is also related to a Gushel--Mukai fourfold in terms of moduli spaces:

\begin{theorem}(Perry--Pertusi--Zhao \cite{PPZ})\label{ppz} Let $Y$ be a very general Gushel--Mukai fourfold. There exist a Mukai vector $v$ and stability conditions $\sigma$ on the Kuznetsov component $Ku(Y)$, such that
there is an isomorphism
  \[ X\cong M_{v,\sigma}(Ku(Y)) \]
  between a double EPW sextic $X$ and the moduli space $M_{v,\sigma}(Ku(Y))$. The very general double EPW sextic is attained in this way.
  
  (That is, there exists $\MM^\circ_{GM4}\subset\MM_{GM4}$, intersection of a countaby infinite number of dense open subsets, such that over $\MM^\circ_{GM4}$ the double EPW sextic has a modular interpretation.)
  
  \end{theorem}    

\begin{proof} For a general $Y$, the moduli space  $M_{v,\sigma}(Ku(Y))$ exists (and is a smooth projective hyperk\"ahler variety) thanks to
\cite[Theorem 1.7]{PPZ}. For the very general $Y$, the moduli space is identified with a double EPW sextic in \cite[Proposition 5.17]{PPZ}. 
   
For later use, we remark that the construction of \cite{PPZ} is actually done in families over the moduli space of Gushel--Mukai fourfolds.
\end{proof}

\begin{remark} The modular construction induces a morphism
   \[ \nu_{PPZ}\colon\ \ \MM^\circ_{GM4}\ \to\ \MM_{dEPW} \ .\]
   One would expect that this morphism $\nu_{PPZ}$ agrees with the rational map $\nu$ of Theorem \ref{im}. This is not proven in \cite{PPZ}; what is known \cite[Proposition 5.17]{PPZ} is that either the two maps agree, or the modular construction gives the double EPW sextic that is dual to the Iliev--Manivel construction. (For completeness, we should mention the recent work \cite{GLZ}, where it is shown that $\nu_{PPZ}$ actually agrees with $\nu$.)
Either way, the morphism $\nu_{PPZ}\colon\MM^\circ_{GM4} \to \MM_{dEPW}$
   is dominant.
   \end{remark}

As a consequence of the modular construction, we have a good understanding of the (Chow groups and) Chow motive of the very general double EPW sextic:

\begin{theorem}\label{bul} Let $Y$ be a very general Gushel--Mukai fourfold, and let $X$ be the double EPW sextic associated to $Y$ via Theorem \ref{ppz}. There is an inclusion of Chow motives
  \[  h(X)\ \hookrightarrow\ \bigoplus h(Y^2)(\ast)\ \ \ \hbox{in}\ \MM_{\rm rat}\ .\]
 \end{theorem}
 
 \begin{proof} This is essentially B\"ulles' result \cite{Bul}, with the improved bound on the exponent obtained in \cite[Theorem 1.1]{FLV3}. For later use, we remark that the argument proving Theorem \ref{bul} can be done in families relative to a base; this is explained in (and crucial to the argument of) \cite{FLV3}.
 \end{proof}

   \subsection{Generically defined cycles and the Franchettina property}  
   \label{ss:gen}
   
   \begin{definition}\label{frank} Let $\YY\to S$ be a smooth projective morphism, where $\YY, S$ are smooth quasi-projective varieties. We say that $\YY\to S$ has the {\em Franchetta property in codimension $j$\/} if the following holds: for every $\Gamma\in A^j(\YY)$ such that the restriction $\Gamma\vert_{Y_s}$ is homologically trivial for all $s\in S$, the restriction $\Gamma\vert_s$ is zero in $A^j(Y_s)$ for all $s\in S$.
 
 We say that $\YY\to S$ has the {\em Franchetta property\/} if $\YY\to S$ has the Franchetta property in codimension $j$ for all $j$.
 \end{definition}
 
 This property is studied in \cite{BL}, \cite{FLV}, \cite{FLV3}, \cite{FLV2}. The following handy ``spread lemma'' implies that in Definition \ref{frank}, it actually suffices to consider the very general fiber:
 
 \begin{lemma}\label{spread} Let $\YY\to S$ be a smooth projective morphism, and $\Gamma\in A^j(\YY)$. Assume that the restriction $\Gamma\vert_{Y_s}$ is 0 in $A^j(Y_s)$ for the very general $s\in S$. Then the restriction $\Gamma\vert_{Y_s}$ is 0 in $A^j(Y_s)$ for all $s\in S$. 
 \end{lemma}
 
 \begin{proof} This is \cite[Lemma 3.1]{Vo}. A proof can be found in \cite[Proposition 2.4]{V3}.
  \end{proof}

 We can also consider a weaker variant of the Franchetta property, by replacing Chow groups $A^\ast()$ by groups $B^\ast()$ of algebraic cycles modulo algebraic equivalence:
 
  \begin{definition}\label{frank2} Let $\YY\to S$ be a smooth projective morphism, where $\YY, S$ are smooth quasi-projective varieties. We say that $\YY\to S$ has the {\em Franchettina property in codimension $j$\/} if the following holds: for every $\Gamma\in B^j(\YY)$ such that the restriction $\Gamma\vert_{Y_s}$ is homologically trivial for all $s\in S$, the restriction $\Gamma\vert_s$ is zero in $B^j(Y_s)$ for all $s\in S$.
 
 We say that $\YY\to S$ has the {\em Franchettina property\/} if $\YY\to S$ has the Franchettina property in codimension $j$ for all $j$.
 \end{definition}
 
 The Franchettina property is studied in \cite{BoLa}.

 \begin{notation}\label{def:gd} Given a family $\YY\to S$ as above, with $Y:=Y_s$ a fiber, we write
   \[ \begin{split} GDA^j_S(Y):=\ima\Bigl( 
  A^j(\YY)\to A^j(Y)\Bigr)\ ,\\
      GDB^j_S(Y):=\ima\Bigl( 
  B^j(\YY)\to B^j(Y)\Bigr)\ \\   
  \end{split}\]
   for the subgroups of {\em generically defined cycles}. 
  In a context where it is clear to which family we are referring, the index $S$ will often be suppressed from the notation.
  \end{notation}
  
  With this notation, the Franchetta (and Franchettina) property amounts to saying that $GDA^\ast_S(Y)$ (resp. $GDB^\ast_S(Y)$) injects into cohomology, under the cycle class map. 
  
  The following is the main result of this subsection:
   
   \begin{proposition}\label{franche} Let $\XX\to\MM^\circ_{GM4}$ be the universal family of double EPW sextics, where $\MM^\circ_{GM4}$ is as in Theorem \ref{ppz}. This family has the Franchettina property, i.e.
   the cycle class map induces injections
       \[ GDB^\ast_{\MM^\circ_{GM4}}(X)\ \to\ H^\ast(X,\QQ) \ ,\]
       for any fiber $X$.   
   \end{proposition}
   
   \begin{proof} By assumption, $X$ is constructed as a moduli space as in Theorem \ref{ppz}. Let us write $\YY\to\MM^\circ_{GM4}$ for the universal family of 
   Gushel--Mukai fourfolds restricted to $\MM^\circ_{GM4}$.
   The modular construction (Theorem \ref{ppz}) and the relation of Chow motives (Theorem \ref{bul}) exist on the level of families over $\MM^\circ_{GM4}$, and so there is a commutative diagram
   \[  \begin{array}[c]{ccc}
                    B^\ast(\XX) & \to & \bigoplus B^\ast(\YY\times_{\MM^\circ_{GM4}}\YY)\\
                    &&\\
                    \downarrow&&\downarrow\\
                    &&\\
                    B^\ast(X) & \hookrightarrow & \bigoplus B^\ast(Y^2)\\
                    &&\\
                    \downarrow&&\downarrow\\
                    &&\\
                    H^\ast(X,\QQ)& \to&\bigoplus H^\ast(Y^2,\QQ)\ ,\\
                         \end{array}\]
                         where the middle horizontal arrow is injective thanks to Theorem \ref{bul}.
                         
                         We are thus reduced to proving the Franchettina property for $Y^2$. This is readily done; the argument is very similar to that for Gushel--Mukai sixfolds in \cite{BoLa}.
   
% First, let us show the Franchetta property for $Y$. 
 Let us denote by
   \[ \BB\subset \bar{\BB}:= \PP H^0\bigl(\Gr(2,5), \OO(1)\oplus \OO(2)\bigr) \]
   the dense open parametrizing smooth (ordinary) Gushel--Mukai fourfolds, and let $\bar{\YY}\to\bar{\BB}$ denote the universal family of all (possibly singular and not dimensionally transverse) complete intersections.
   As there is a surjection $\BB\to\MM_{GM4}$, it will suffice to establish the Franchettina property for $\YY\times_{\bar{\BB}} \YY\to B$. The fact that the line bundles $\OO(1)$ and $\OO(2)$ are very ample on $\Gr(2,5)$ implies that the morphism
   \[    \bar{\YY}\times_{\bar{\BB}} \bar{\YY}\ \to\  \Gr(2,5)\times\Gr(2,5) \]
   is a {\em stratified projective bundle\/} (in the sense of \cite{FLV3}), with strata the diagonal $\Delta_{\Gr(2,5)}$ and its complement.
   Then the stratified projective bundle argument \cite[Proposition 2.6]{FLV3} gives an equality
    \begin{equation}\label{GDAYY} GDA^\ast_{\BB}(Y\times Y) =  \Bigl\langle  (p_i)^\ast  GDA^\ast_{\BB}(Y), \Delta_Y \Bigr\rangle\ ,\end{equation}
    where $p_i\colon Y\times Y\to Y$, $i=1,2$, denote the two projections.
    
    In order to better understand the right-hand side of equality \eqref{GDAYY}, we establish 3 lemmata:
    
    \begin{lemma}\label{l1} Let $Y$ be a smooth, ordinary Gushel--Mukai fourfold. One has
      \[  GDA^\ast_{\BB}(Y) =\bigl\langle h,c\bigr\rangle\ ,\]
      where $h\in A^1(Y)$ is the Pl\"ucker hyperplane section and $c\in A^2(Y)$ is the restriction of $c_2(Q)$ with $Q$ the tautological bundle on $\Gr(2,5)$.
      
      Moreover, $GDB^\ast_{\BB}(Y) $ injects into cohomology.    
        \end{lemma}
      
      \begin{proof} The morphism $\bar{\YY}\to\Gr(2,5)$ is a projective bundle. The projective bundle formula, plus the argument found in \cite[Proof of Lemma 1.1]{PSY}, implies
        \[   GDA^\ast_{\BB}(Y) =\ima\Bigl( A^\ast(\Gr(2,5))     \to A^\ast(Y)\Bigr)\ .\]
        Since the Chow ring of $\Gr(2,5)$ is generated by the Pl\"ucker hyperplane and $c_2$  of the tautological bundle \cite{3264}, this gives the expression for the generically defined cycles.
%        For the injectivity, as the Picard number of $Y$ is 1 we have a cohomological relation 
%          \[ c\cdot h = \beta h^3\ \ \ \hbox{in}\ H^{6}(Y,\QQ)\ ,\]
%          for some $\beta\in\QQ$.
%      To lift this to rational equivalence,  let us write $Y=Y^\prime\cap Q$, where $Q$ is a quadric and $Y^\prime=\Gr(2,5)\cap H$ is the intersection of the Grassmannian with a Pl\"ucker hyperplane section. Let $j\colon Y\to Y^\prime$ denote the inclusion morphism. We have
%          \[   c\cdot h
     
     For the injectivity, this is true for trivial reasons: $Y$ is a Fano fourfold and hence $A^4_{hom}(Y)=0$, and so (as can be se seen from the Bloch--Srinivas argument \cite{BS}) we have
     $B^\ast_{hom}(Y)=0$.
         \end{proof}

        \begin{lemma}\label{l2}  Let $Y$ be a smooth, ordinary Gushel--Mukai fourfold. One has equality
       \[   \Delta_Y\cdot (p_i)^\ast(h) =  p\ \ \ \hbox{in}\ A^5(Y\times Y)   \ ,\]
       where $p$ is a polynomial in $(p_i)^\ast(h)$ and $(p_i)^\ast(c)$.
              \end{lemma}
              
   \begin{proof} Let us write $Y=Y^\prime\cap Q$, where $Q$ is a quadric and $Y^\prime=\gr(2,5)\cap H$ is the intersection of the Grassmannian with a Pl\"ucker hyperplane section. The excess intersection formula \cite[Theorem 6.3]{F} gives an equality
     \[    \Delta_Y\cdot (p_i)^\ast(h) = {1\over 2}\ \Delta_{Y^\prime}\vert_Y\ \ \ \hbox{in}\ A^5(Y\times Y)\ .\]
     But the fivefold $Y^\prime$ has trivial Chow groups \cite[Proposition 4.6]{GrPf}, and so (by weak Lefschetz, applied to the inclusion $Y^\prime\hookrightarrow\Gr(2,5)$) we know that $\Delta_{Y^\prime}$ is a polynomial in the classes $h\in A^1(\Gr(2,5))$, $c\in A^2(\Gr(2,5))$ restricted to $Y^\prime$ and pulled-back to $Y^\prime\times Y^\prime$. This settles the lemma.
              \end{proof}

\begin{lemma}\label{l3}  Let $Y$ be a smooth, ordinary Gushel--Mukai fourfold. One has equality
          \[   \Delta_Y\cdot (p_i)^\ast(c) =  q\ \ \ \hbox{in}\ B^6(Y\times Y)   \ ,\]
       where $q$ is a polynomial in $(p_i)^\ast(h)$ and $(p_i)^\ast(c)$.         
        \end{lemma}       
       
\begin{proof} First, let us remark that one has a relation on the level of cohomology
  \[  \Delta_Y\cdot (p_i)^\ast(c) =  q\ \ \ \hbox{in}\ H^{12}(Y\times Y,\QQ) \]      
  (where $q$ is as in the lemma). Indeed, the correspondence $\Delta_Y\cdot (p_i)^\ast(c) $
 acts on $H^\ast(Y,\QQ)$ as cupping with the class $c$. As the primitive cohomology of $X$ is concentrated in degree $4$, the correspondence $\Delta_Y\cdot (p_i)^\ast(c) $ acts as 0 on the primitive cohomology. The algebraic part of $H^\ast(Y,\QQ)$ can be expressed in terms of $h$ and $c$, and so there exists a polynomial $q$ in $h$ and $c$ such that $q$ and  $\Delta_Y\cdot (p_i)^\ast(c) $ act in the same way on $H^\ast(Y,\QQ)$. Manin's identity principle (plus the K\"unneth formula in cohomology) then implies that
   \[    \Delta_Y\cdot (p_i)^\ast (c) = q\ \ \hbox{in}\ H^{12}(Y\times Y,\QQ)\ .\] 
   
It remains to upgrade to algebraic equivalence. Since $Y$ is a Fano fourfold, we have $A^4_{hom}(Y)=0$. The Bloch--Srinivas argument (in the precise form of \cite[Theorem 3.11]{Vab}) then gives us a surface $S$ and an inclusion of Chow motives
  \[  h(Y)\ \hookrightarrow\ h(S)(-1)\oplus \bigoplus\one(\ast)\ \ \ \hbox{in}\ \MM_{\rm rat}.\]
  This induces in particular an inclusion
  \[ B^6_{hom}(Y \times Y) \ \hookrightarrow \  B^4_{hom}(S\times S)\ .\]
But $B^4_{hom}(S\times S)=0$ (homological and algebraic equivalence coincide for zero-cycles), and so $B^6_{hom}(X\times X)=0$.  
The lemma is proven.
   \end{proof}
       
       With the aid of the last two lemmata, the equality \eqref{GDAYY} reduces to
       \[   GDB^\ast_{\BB}(Y\times Y) =  \Bigl\langle  (p_i)^\ast  GDB^\ast_{\BB}(Y)\Bigr\rangle \oplus \QQ[ \Delta_Y]  \ .\]
      Now we observe that the class of the diagonal $\Delta_Y$ in cohomology is linearly independent from the decomposable classes $ \bigl\langle  (p_i)^\ast  GDB^\ast_{\BB}(Y)\bigr\rangle  $ (indeed, the diagonal acts as the identity on $H^{3,1}(Y)$, whereas the decomposable classes act as zero on $H^{3,1}(Y)$).
      Using the K\"unneth formula in cohomology, one is thus reduced to proving the Franchettina property for $\YY\to\BB$; this is Lemma \ref{l1}.          
               \end{proof}
          
  \begin{remark}\label{pity} The reason for working with algebraic equivalence (instead of rational equivalence) is that I have not been able to prove the injectivity statement of Lemma \ref{l1} and Lemma \ref{l3} modulo rational equivalence. That is, if one is able to upgrade these two statements to rational equivalence, one obtains the Franchetta property (rather than the Franchettina property) for 
  double EPW sextics.
   \end{remark}      
  
  We will make use of the Franchettina result (Proposition \ref{franche}) in the following way:  
  
  \begin{corollary}\label{Fr2} Let $Y$ be a Gushel--Mukai fourfold in $\MM^\circ_{GM4}$, and let $X$ be the double EPW sextic associated to $Y$ via Theorem \ref{ppz}. We have
    \[ GDB^2_{\MM^\circ_{GM4}}(X) =  \QQ[h^2]\oplus \QQ[c]\ .\]
  (Here, $h$ and $c$ are as in Theorem \ref{fer}.)
      \end{corollary}
  
  \begin{proof} 
     The fourfold $X$ is a deformation of a Hilbert square of a K3 surface, and so cup product induces an isomorphism
    \[ \hbox{Sym}^2 H^2(X,\QQ)\ \xrightarrow{\cong}\ H^4(X,\QQ)\ .\]
    It follows that
    \[ H^{2,2}(X,\QQ)= A^1(X)\cdot A^1(X)\oplus \QQ[c]\ ,\]
    and so for the very general $X$ the space of Hodge classes $H^{2,2}(X,\QQ)$ is two-dimensional. Since $ GDB^2_{\MM^\circ_{GM4}}(X)$ injects into cohomology (Proposition \ref{franche}), this proves (\rom1). 
  \end{proof}

  \subsection{A quadratic relation}
  
  \begin{definition} Given a general Gushel--Mukai fourfold $Y$, let $F=F(Y)$ be the variety of conics. We define 
      \[I_F:= \bigl\{ (c,c^\prime)\in F\times F\ \big\vert \ c\cap c^\prime\not=\emptyset\bigr\}\ \ \  \subset F\times F\] 
      as the subvariety of pairs of intersecting conics.
      
  We also define 
      \[  W:= \bigl\{  (c,c^\prime)   \ \big\vert \exists V_4 \hbox{\ such\ that\  $c$\ and\ $c^\prime$\ are\ residual\ in\ $S_{V_4}$}  \bigr\}\ \ \ \subset F\times F\ ,\]
 where $V_4\subset V_5$ is a 4-dimensional subvector space, and $S_{V_4}$ is the surface $\Gr(2,V_4)\cap H\cap Q$ contained in $Y:=\Gr(2,V_5)\cap H\cap Q$.
   \end{definition}
   
   \begin{remark} The variety $I_F\subset F\times F$ is 8-dimensional. There is equality
      \begin{equation}\label{Ipq}  I_F = (p\times p)_\ast (q\times q)^\ast (\Delta_Y) = {}^t P\circ P\ \ \ \hbox{in}\ A^2(F\times F)\ .\end{equation}
      The variety $W$ is 6-dimensional.
              \end{remark}
              
   \begin{lemma}\label{w} Let $Y$ be a general Gushel--Mukai fourfold, with $F=F(Y)$ the variety of conics and $\pi\colon F\to X$ the morphism to the associated double EPW sextic.
   Let $h_F\in A^1(F)$ be an ample divisor, and let $\iota\in\aut(X)$ denote the covering involution.
   One has equality
     \[ (\pi\times\pi)_\ast \bigl(W\cdot (h_F\times h_F)\bigr)   = d\, \Gamma_\iota + R\ \ \ \hbox{in}\ A^4(X\times X)\ ,\]
     where $d\in\QQ^\ast$ is some constant depending on $h_F$, and $R$ is supported on $x_\rho\times X\cup x_\sigma\times X\cup X\times x_\rho\cup X\times x_\sigma$.
        \end{lemma}  
   
   \begin{proof} This follows from the explicit description of $\iota$ in terms of conics given in \cite[Section 3.4]{Zh}: if $x\in X$ is a general point, a conic $c\in\pi^{-1}(x)$ will be a $\tau$-conic. The conic $c$ is contained in a unique surface $S_{V_4}$ inside the Gushel--Mukai fourfold $Y$. Let $c^\prime$ be the residual conic in $S_{V_4}$ (i.e. $S_{V_4}\cap H=c\cup c^\prime$). Then
     \[ \iota(x)=\pi(c^\prime)\ .\]
    \end{proof}

  \begin{proposition}\label{zh}(Zhang \cite{Zh}) Let $Y$ be a general Gushel--Mukai fourfold, and let $F=F(Y)$ be the variety of conics contained in $Y$. There is a relation
    \[  (I_F)^2= \alpha W + I_F\cdot A + B + C\ \ \ \hbox{in}\ A^4(F\times F)\ ,\]
    where $\alpha\in\QQ^\ast$, $A$ and $B$ are decomposable and generically defined correspondences, i.e.
      \[ A,B\ \ \in\ \Bigl\langle   (p_i)^\ast  GDA^\ast_{\MM_{GM4}}(F)\Bigr\rangle\ \] 
      (with $p_i$ the projection from $F\times F$ to the $i$-th factor), and $C$ is supported on $\Sigma_2:=(F_\rho(Y)\cup F_\sigma(Y))\times (F_\rho(Y)\cup F_\sigma(Y))$.
    \end{proposition}
    
    \begin{proof} This is essentially \cite[Proposition 4.10]{Zh}, whose argument is modelled on the analogous statement for the Beauville--Donagi fourfolds proven by Voisin \cite[Proposition 3.3]{V17}. 
    The fact that $A$ and $B$ are generically defined is not explicitly stated in loc. cit., but follows from inspection of the proof: indeed, the proof is a Chern class computation, and $A$ and $B$ are
    expressions in Chern classes of the tangent bundle of $Y$ and Chern classes of the relative tangent bundle of $P\to F$, which are all generically defined (with respect to $\MM_{GM4}$).
     \end{proof}
    
    For our purposes, it will be more convenient to obtain a quadratic relation on $X\times X$ rather than on $F\times F$ (one reason being that the variety of conics $F$ is well-behaved only for the {\em general\/} double EPW sextic; another reason being that we want to involve the covering involution $\iota$, which exists on $X$ but not (a priori) on $F$). To this end, we propose the following definition:
    
   \begin{definition} Let $Y$ be a Gushel--Mukai fourfold that is general in the sense of Theorem \ref{im}, let $F=F(Y)$ be the variety of conics, and let $\pi\colon F\to X$ be the morphism to the associated
   double EPW sextic. We define the subvariety $I_\pi\subset X\times X$ by
     \[  I_\pi:= (\pi\times\pi)(I_F)\ \ \ \subset\ X\times X\ .\]
      \end{definition} 
      
      \begin{remark} We note that a priori, the subvariety $I_\pi$ is {\em not \/} intrinsic to $X$. That is, if $Y_1, Y_2$ are two Gushel--Mukai fourfolds with the same associated double EPW sextic $X$
      (such $Y_1$ and $Y_2$ are called ``period partners''), we do not know whether $I_{\pi_1}$ and $I_{\pi_2}$ coincide.
      
       We also remark that if $(c,c^\prime)\in I_F$ are intersecting conics that are both $\tau$-conics, then $I_F$ contains the whole fiber $(\pi\times\pi)^{-1}(\pi(c),\pi(c^\prime))$.
      Thus, we have a relation
      \begin{equation}\label{fiber}   I_F = (\pi\times\pi)^\ast (I_\pi) + R\ \ \ \hbox{in}\ A^2(F\times F)\ ,\end{equation}
      where $R$ is supported on $\bigl((F_\rho(Y)\cup F_\sigma(Y))\times F\bigr) \cup \bigl(F\times (F_\rho(Y)\cup F_\sigma(Y))\bigr)$.      
        \end{remark}

    \begin{proposition}\label{zh2} Let 
    $Y$ be a Gushel--Mukai fourfold that is general in the sense of Theorem \ref{im}, let $X$ be the associated double EPW sextic, and let $\iota\in \aut(X)$ denote the covering involution. There is a relation
      \[ (I_\pi)^2 = \alpha \Gamma_\iota + I_\pi\cdot A + B\ \ \ \hbox{in}\ A^4(X\times X)\ ,\]
      where $\alpha\in\QQ^\ast$, and $A$ and $B$ 
      are as in Proposition \ref{zh}.
     \end{proposition}
    
    \begin{proof} Let $h_F\in A^1(F)$ be an ample divisor. Relation \eqref{fiber} implies 
      \begin{equation}\label{d} I_\pi= {d}\, (\pi\times\pi)_\ast \bigl(I_F\cdot (h_F\times h_F)\bigr) \ \ \ \hbox{in}\ A^2(X\times X)\ ,\end{equation}
      where $d\in\QQ^\ast$ depends on $h_F$.
      
      Intersecting the quadratic relation on $F\times F$ (Proposition \ref{zh}) with $h_F\times h_F$ and pushing forward to $X\times X$, we obtain
      \[   (\pi\times\pi)_\ast \Bigl( I_F\cdot (h_F\times h_F) \cdot ((\pi\times\pi)^\ast (I_\pi) + R)  \Bigr) =     (\pi\times\pi)_\ast \Bigl(  \alpha W + I_F\cdot A + B + C\Bigr)\ \ \ \hbox{in}\ A^4(X\times X)\ .\]
      Applying the projection formula and relation \eqref{d} to the left-hand side, we obtain
       \[ {d}\, (I_\pi)^2 +R^\prime =    (\pi\times\pi)_\ast \Bigl(  (\alpha W + I_F\cdot A + B + C)\cdot (h_F\times h_F)\Bigr)\ \ \ \hbox{in}\ A^4(X\times X)\ ,\]
       where $R^\prime$ is supported on $x_\rho\times X\cup x_\sigma\times X\cup X\times x_\rho\cup X\times x_\sigma$. This implies that $R^\prime$ is decomposable and generically defined (with respect to $\MM_{GM4}$). We observe that
       \[  (\pi\times\pi)_\ast \Bigl( (B+C)\cdot (h_F\times h_F)\Bigr) \]
       is also decomposable and generically defined (with respect to $\MM_{GM4}$). Hence, we can rewrite the above equality as
        \[ { d}\, (I_\pi)^2 =    (\pi\times\pi)_\ast \Bigl(  (\alpha W + I_F\cdot A)\cdot (h_F\times h_F)\Bigr)    + B^\prime \ \ \ \hbox{in}\ A^4(X\times X)\ ,\]
        for some $B^\prime$ that is decomposable and generically defined. Again using relation \eqref{fiber} and the projection formula, we find that
        \[   (\pi\times\pi)_\ast \Bigl(  I_F\cdot A\cdot (h_F\times h_F)\Bigr) =   {1\over d}\, I_\pi\cdot      (\pi\times\pi)_\ast (A) + R^{\prime\prime} \ ,\]
        where $R^{\prime\prime}$ is again (supported on $x_\rho\times X\cup x_\sigma\times X\cup X\times x_\rho\cup X\times x_\sigma$ and hence) decomposable and generically defined. The above equality thus simplifies to
        \[       (I_\pi)^2 =  \alpha^\prime  (\pi\times\pi)_\ast \Bigl(  W \cdot (h_F\times h_F)\Bigr)  + I_\pi\cdot A^\prime  + B^{\prime\prime} \ \ \ \hbox{in}\ A^4(X\times X)\ .\]
        
       Finally, applying Lemma \ref{w} the equality further simplifies to
                \[       (I_\pi)^2 =  \alpha^{\prime\prime}  \Gamma_\iota  + I_\pi\cdot A^\prime  + B^{\prime\prime\prime} \ \ \ \hbox{in}\ A^4(X\times X)\ ,\]
                with $ \alpha^{\prime\prime}\in\QQ^\ast$ and $A^\prime,B^{\prime\prime\prime}$ correspondences of the required type.
          \end{proof}

  \subsection{A linear relation}
  
  \begin{proposition}\label{I} Let $X$ be a double EPW sextic that is general in the sense of Theorem \ref{im}, and let $\iota$ denote the covering involution. There is a relation
    \[ I_\pi + (\iota,\ide)^\ast (I_\pi) = Z\ \ \ \hbox{in}\ A^2(X\times X)\ ,\]
    where $Z$ is decomposable and generically defined, that is 
      \[  Z\ \ \in\ \Bigl\langle (p_1)^\ast GDA_{\MM_{GM4}}^\ast(X), (p_2)^\ast GDA_{\MM_{GM4}}^\ast(X)\Bigr\rangle\ .\]
   \end{proposition}
   
   \begin{proof} This is modelled on \cite[Lemma 3.12]{V17}, which is the analogous result for the Beauville--Donagi fourfolds.
   
 By hypothesis, $X$ can be constructed via the Iliev--Manivel construction (Theorem \ref{im}) starting from a Gushel--Mukai fourfold $Y$. Given a general point $x\in X$, let us write $C_x\in A^3(Y)$ for the class of the linear system of conics on $Y$ corresponding to $x$. We have an equality
    \[  C_x + C_{\iota(x)}=  C\ \ \ \hbox{in}\ A^3(Y)\ ,\]
    where $C\in A^3(Y)$ is some class independent of the point $x$
    \cite[Lemma 4.3]{Zh}.
    Applying the universal conic $P$ as a correspondence, it follows that there is also an equality
    \[  p_\ast q^\ast (C_x) + p_\ast q^\ast (C_{\iota(x)}) = p_\ast q^\ast (C)\ \ \ \hbox{in} \ A^2(F)\ .\]
   In view of equality \eqref{Ipq}, this can be interpreted as saying
     \[ I_F\vert_{y\times F} + I_F\vert_{ y^\prime\times F} =     p_\ast q^\ast (C)\ \ \ \hbox{in} \ A^2(F)\ ,\]   
     where $y$ and $y^\prime\in F$ are points mapping to $x$ resp. to $\iota(x)$. Intersecting with $F\times h_F$ (where $h_F\in A^1(F)$ is ample) and pushing forward to $X\times X$, it follows that there is equality
     \[ I_\pi\vert_{x\times X} +  \bigl((\iota,\ide)^\ast( I_\pi)\bigr)\vert_{\iota(x)\times X} = C^\prime\ \ \ \hbox{in}\ A^2(X) \ ,\]
     for some constant $C^\prime\in A^2(X)$ independent of $x$. As this applies to a general point $x$, it follows that there is an equality of the action of correspondences
     \[  \bigl(I_\pi +  (\iota,\ide)^\ast (I_\pi)\bigr){}_\ast   = C^\prime_\ast\colon\ \ \ A^4(X)\ \to\ A^2(X)\ .\]
     Then, the ``decomposition of the diagonal'' argument of Bloch--Srinivas \cite{BS} implies that there is equality of correspondences
     \[  I_\pi + (\iota,\ide)^\ast (I_\pi) = C^\prime\times X + \Gamma\ \ \ \hbox{in}\ A^2(X\times X)\ \]
     for some $\Gamma$ supported on $X\times D$, where $D\subset X$ is a divisor. Let $\wt{D}\to D$ be a resolution of singularities. Since $\pic^0(X)=0$, any divisor on $X\times\wt{D}$ is a linear combination of pullbacks of divisors on $X$ and on $\wt{D}$. This implies that the cycle $\Gamma$ is decomposable, and so we obtain 
     \begin{equation}\label{z}  I_\pi + (\iota,\ide)^\ast (I_\pi) = Z\ \ \ \hbox{in}\ A^2(X\times X)\ ,\end{equation}
     with $Z$ decomposable.
     
   It remains to see that $Z$ is generically defined. To this end, we remark that the left-hand side $I_\pi + (\iota,\ide)^\ast (I_\pi)$ is obviously generically defined (with respect to $\MM_{dEPW}$ and so a fortiori with respect to $\MM_{GM4}$). We know that $Z$ is of the form
     \[ Z= C_1\times X + \sum_i D_i\times D_i^\prime + X\times C_2\ \ \ \hbox{in}\ A^2(X\times X)\ ,\]
  with $C_i\in A^2(X)$ and $D_i, D_i^\prime\in A^1(X)$. By intersecting both sides of the relation \eqref{z}
 %  \[ I_\pi + (\iota,\ide)^\ast (I_\pi) = Z\ \ \ \hbox{in}\ A^2(X\times X)\ \]
   with $X\times h^4$ and pushing forward under $p_1$, we find that $C_1$ is generically defined. Likewise, we find that $C_2$ is generically defined. Restricting to the very general $X$, we may assume $A^1(X)=\QQ$ and so all the $D_i$ and $D_i^\prime$ are proportional to $h$. That is, we have now proven the proposition for the very general $X$.

To extend to {\em all\/} double EPW sextics attained by the Iliev--Manivel construction, one can apply the spread lemma (Lemma \ref{spread}) to conclude.
     \end{proof}

\section{Proof of the main result}

This section contains the proof of the main result, as announced in the introduction:

\begin{theorem}\label{main} Let $X$ be a smooth double EPW sextic. The $\QQ$-subalgebra
  \[  \bigl\langle A^1(X), c_j(X)\bigr\rangle\ \ \subset\ A^\ast(X) \]
  injects into cohomology, under the cycle class map.
  
  Moreover, let $A^2(X)^+\subset A^2(X)$ denote the subgroup of cycles invariant under the covering involution $\iota$. The cycle class map induces injections
  \[  \bigl\langle A^1(X), c_j(X), A^2(X)^+\bigr\rangle\ \cap A^i(X)\ \hookrightarrow\ H^{2i}(X,\QQ) \ \ \ \hbox{for}\ i\ge 3\ .\]
  \end{theorem}
  
  \begin{proof} In view of Theorem \ref{lv}, it only remains to prove that the cycle class map induces an injection
   \[  \bigl\langle A^1(X), c_j(X), A^2(X)^+\bigr\rangle\ \cap A^3(X)\ \hookrightarrow\ H^{6}(X,\QQ)\ .\]
   What's more, as noted in \cite{LV}, it suffices to prove that there is an inclusion
   %there exists $\beta_D\in\QQ$ such that
   %\begin{equation}\label{need1}  D^3 =  \beta_D h^2\cdot D\ \ \hbox{in}\ A^3(X)\ ,\end{equation}
   %and 
   \begin{equation}\label{need}  A^1(X)^-\cdot A^2(X)^+ \ \subset\  A^1(X)^-\cdot h^2\ .\end{equation}
   Because of the hard Lefschetz theorem, the wanted relation \eqref{need} holds true in cohomology, and so it remains to prove that
   \begin{equation}\label{need2}    A^1(X)^-\cdot A^2(X)^+\cap A^3_{hom}(X)^- = 0\ .\end{equation}
   
   What follows is closely modelled on the argument for codimension 3 cycles on the Beauville--Donagi fourfold developed by Voisin (cf. \cite[Section 3]{V17}).     
   At the heart of the argument is the following proposition:
   
   \begin{proposition}\label{heart} Let $X$ be a smooth double EPW sextic, and let $\iota$ denote its covering involution.
   There is a relation
     \[ {1\over 2}\Bigl(\Delta_X -\Gamma_\iota\Bigr)= J_1+\cdots + J_k + B\ \ \ \hbox{in}\ A^4(X\times X)\ ,\]
   with the following properties:
   
   \begin{itemize}
   \item $B$ is decomposable and generically defined, i.e. 
       \[B\ \ \in\ \bigl\langle (p_1)^\ast GDA^\ast_{\MM_{dEPW}}(X),  (p_2)^\ast GDA^\ast_{\MM_{dEPW}}(X)     \bigr\rangle\ ;\]
    \item
    each $J_i$ is a composition of the correspondences  $   {1\over 2}\Bigl(\Delta_X -\Gamma_\iota\Bigr)$ and $I\cdot A$, where \[ \begin{split} &I\in GDA^2_{\MM_{dEPW}}(X\times X)\ ,\\
    &A\ \ \in\ \Bigl\langle (p_1)^\ast GDA^\ast_{\MM_{dEPW}}(X),  (p_2)^\ast GDA^\ast_{\MM_{dEPW}}(X)     \Bigr\rangle\cap A^2(X\times X)\ ,\\
         \end{split}\]
     with $I\cdot A$ occurring at least once in each $J_i$.
      \end{itemize}
        \end{proposition}
   
   Before embarking on the proof of Proposition \ref{heart}, let us first show that Proposition \ref{heart} implies the wanted equality \eqref{need2} (and hence the theorem). 
%   The celebrated nilpotence theorem of Voisin \cite{V2} and Voevodsky \cite{Voe} turns the relation of Proposition \ref{heart} into a relation modulo rational equivalence
%   \begin{equation}\label{eq4}  {1\over 2}\Bigl(\Delta_X -\Gamma_\iota\Bigr)= J + B\ \ \ \hbox{in}\ A^4(X\times X)\ ,\end{equation}
%   where $B$ is decomposable and $J$ is a sum of compositions of correspondences of the form
%   \[ J=    {1\over 2}\Bigl(\Delta_X -\Gamma_\iota\Bigr)\circ (I\cdot A)\circ {1\over 2}\Bigl(\Delta_X -\Gamma_\iota\Bigr)\circ \cdots\circ (I\cdot A)\circ {1\over 2}\Bigl(\Delta_X -\Gamma_\iota\Bigr)\ ,\]   
%         in which $I\cdot A$ occurs at least once.
We consider an element
    \[  \gamma\in A^1(X)^-\cdot A^2(X)^+\cap A^3_{hom}(X)^-\ .\]
    To prove the wanted equality \eqref{need2}, we need to prove that $\gamma$ is 0. Applying the equality of correspondences of Proposition \ref{heart} to the cycle $\gamma$, we obtain that
    \[   \gamma = (J_1 + \cdots + J_k+ B)_\ast(\gamma)\ \ \ \hbox{in}\ A^3(X)\ .\]
    The correspondence $B$, being decomposable, acts as zero on the homologically trivial cycle $\gamma\in A^3_{hom}(X)$, and so this equality further simplifies to   
    \[  \gamma =    (J_1 + \cdots + J_k)_\ast (\gamma) \ \ \ \hbox{in}\ A^3(X)\ .\]  
    Since $  {1\over 2}\Bigl(\Delta_X -\Gamma_\iota\Bigr)$ acts as the identity on $\gamma$, and each $J_i$ contains the correspondence $I\cdot A$, 
    this can be rewritten as
      \[ \gamma =  (\hbox{something})_\ast (I\cdot A)_\ast (\gamma)  \ \ \ \hbox{in}\ A^3(X)\ .\]      
   The desired equality \eqref{need2} now follows from the following lemma:
    
      \begin{lemma}\label{end} Let $I, A$ be as in Proposition \ref{heart}. Then
     \[ (I\cdot A)_\ast(\gamma)  =0\ \ \ \forall\ \gamma\in A^1(X)^-\cdot A^2(X)^+\cap A^3_{hom}(X)^-\ .  \]
     \end{lemma}
     
   \begin{proof} (of the lemma) We need to remember that the correspondence $A\in A^2(X\times X)$ is {\em decomposable\/}, i.e. it is a linear combinations of pullbacks to $X\times X$ of classes in $A^1(X)$ and $A^2(X)$. Hence, we are reduced to proving that correspondences of the form
   %linear combinations of pullbacks to $X\times X$ of the classes $h\in A^1(X)$ and $c\in A^2(X)$. Hence, we are reduced to proving that the correspondences
  % \[  I\cdot (p_1)^\ast(h^2),\  I\cdot (p_1)^\ast(c),\  I\cdot (p_1)^\ast(h)\cdot(p_2)^\ast(h), \ I\cdot (p_2)^\ast(h^2), \ I\cdot (p_2)^\ast(c) \]
    \[  I\cdot (p_1)^\ast(a),\  I\cdot (p_2)^\ast(a),\ I\cdot (p_1)^\ast(D)\cdot (p_2)^\ast(D^\prime) \]
   act as zero on the cycle $\gamma$, where $a\in A^2(X)$ and $D,D^\prime\in A^1(X)$.
  
  Let us now treat these correspondences case by case. For the first, the projection formula gives us equality
  \[ \bigl(    I\cdot (p_1)^\ast(a)\bigr){}_\ast(\gamma) =  I_\ast (\gamma\cdot a)\ ,\]
  which is zero for reasons of dimension. 
   For the second correspondence, the projection formula gives
  \[ \bigl(  I\cdot (p_2)^\ast(a)\bigr){}_\ast(\gamma)= a\cdot I_\ast(\gamma)\ .\]
  But $I_\ast(\gamma)$ is in $A^1(X)$; since it is homologically trivial and $A^1_{hom}(X)=0$ we find that $I_\ast(\gamma)=0$.  Finally, for the third correspondence, the projection formula implies
  \[ \bigl(  I\cdot (p_1)^\ast(D)\cdot(p_2)^\ast(D^\prime)\bigr){}_\ast(\gamma) =  D^\prime\cdot I_\ast(\gamma\cdot D)\ .\]
  We know from Theorem \ref{lv} that $\gamma\cdot D$ is proportional to $h^4$; since it is also homologically trivial we find that $\gamma\cdot D=0$.
  The lemma is proven.
  \end{proof}
    
    It remains to prove Proposition \ref{heart}. Since all cycles in the equality are generically defined (with respect to $\MM_{dEPW}$), in view of the spread lemma (Lemma \ref{spread}) it suffices to prove the equality of Proposition \ref{heart} for a very general double EPW sextic $X$. We may, and will, thus assume that $X$ is attained both by the Iliev--Manivel construction (Theorem \ref{im}) and by the modular construction (Theorem \ref{ppz}), and prove the proposition in this case. A further reduction step is to reduce to {\em algebraic equivalence\/}, and prove the
    equality
    \begin{equation}\label{onlyalg}  {1\over 2}\Bigl(\Delta_X -\Gamma_\iota\Bigr)= I\cdot A + B\ \ \ \hbox{in}\ B^4(X\times X)\ ,\end{equation}   
    where $I$ and $A$ and $B$ are as in Proposition \ref{heart}. The algebraic equivalence \eqref{onlyalg} implies the rational equivalence of Proposition \ref{heart} by applying the celebrated nilpotence theorem of Voevodsky \cite{Voe} and Voisin \cite{V2}: indeed, the nilpotence theorem transforms \eqref{onlyalg} into an equality
    \[  \Bigl(   {1\over 2}(\Delta_X -\Gamma_\iota)- I\cdot A - B    \Bigr)^{\circ M}=0\ \ \ \hbox{in}\ A^4(X\times X) \]
    for some $M\in\NN$. Developing this expression, and observing that $ {1\over 2}(\Delta_X -\Gamma_\iota)$ is idempotent while any generically defined correspondence composed with $B$ is decomposable and generically defined, one obtains the equality of Proposition \ref{heart}.
    
    To prove equality \eqref{onlyalg}, let us start with the relation given by Proposition \ref{I} and reduce to algebraic equivalence:
    \[ (\iota,\ide)^\ast I_\pi = -I_\pi+Z\ \ \ \hbox{in}\ B^2(X\times X)\ .\]
    Taking squares on both sides, this relation implies
    \begin{equation}\label{square}  \begin{split} (\iota,\ide)^\ast( I_\pi^2)   &= \bigl( (\iota,\ide)^\ast I_\pi\bigr)^2   \\
                                                                                                &= I_\pi^2-2I_\pi\cdot Z +Z^2\ \ \ \ \hbox{in}\ B^4(X\times X)\ .\\
                                                                                                \end{split} \end{equation}
       But we already have another relation involving $I_\pi^2$, given by Proposition \ref{zh2}. Plugging in the relation of Proposition \ref{zh2} into the equality \eqref{square}, we obtain
       \[      (\iota,\ide)^\ast \Bigl(     \alpha \Gamma_\iota + I_\pi\cdot A + B\Bigr) =\alpha \Gamma_\iota + I_\pi\cdot A + B -2I_\pi\cdot Z +Z^2\ \ \ \ \hbox{in}\ B^4(X\times X)\ .\]
       Upon developing (and observing that $(\iota,\ide)^\ast (\Gamma_\iota)$ is nothing but the diagonal $\Delta_X$), we see that this is equivalent to
         \begin{equation}\label{rewrite} \alpha (\Delta_X-\Gamma_\iota) = -(\iota,\ide)^\ast (I_\pi\cdot A) + I_\pi\cdot A +B- 2I_\pi\cdot Z +Z^2 \ \ \ \hbox{in}\ B^4(X\times X)\ .\end{equation}       
       Since $Z$ is decomposable and generically defined, so is $Z^2$. As for the correspondence $A$, we know (Corollary \ref{Fr2}) it is built from $h\in B^1(X)$ and $c\in B^2(X)$ and so $A$ is $\iota$-invariant; we can thus write
       \[  (\iota,\ide)^\ast (I_\pi\cdot A) =(\iota,\ide)^\ast (I_\pi)\cdot A\ \ \ \hbox{in}\ B^4(X\times X)\ .\]
     In view of these two observations, up to changing $A$ and $B$ (but keeping the property that they are generically defined and decomposable), and replacing $I_\pi$ by some other generically defined correspondence $I$, we can rewrite \eqref{rewrite}
     as an equality
     \begin{equation}\label{this} \alpha (\Delta_X-\Gamma_\iota) =  I\cdot A +B \ \ \ \hbox{in}\ B^4(X\times X)\ . \end{equation}
  The only thing left to worry about is the fact that $A$ and $B$ and $I$ in \eqref{this} are generically defined with respect to $\MM_{GM4}$, whereas in \eqref{onlyalg} we require them to be generically defined with respect to $\MM_{dEPW}$. To remedy this, let us consider the base change
    \[  \nu\colon\ \ \MM_{GM4}\ \dashrightarrow \ \MM_{dEPW}\ \]
    given by Theorem \ref{im}.
   Cutting with hyperplane sections and shrinking $\MM_{dEPW}$, we may assume $\nu$ is a finite morphism of degree $d$. That is, we have isomorphic copies $X_1,\ldots,X_d$ of $X$ (corresponding to different Gushel--Mukai fourfolds with the same associated double EPW sextic $X$), and the base change induces a morphism
    \[ \nu\colon\ \ \cup_{i=1}^d X_i\ \to\ X\ .\]
    The above equality \eqref{this} applies to each $X_i$, i.e. we have
    \[    \alpha_i (\Delta_{X_i}-\Gamma_{\iota_i}) =  I_i\cdot A_i +B_i \ \ \ \hbox{in}\ B^4(X_i\times X_i)    \ \ \ (i=1,\ldots,d)\ .   \]
    We may assume the correspondence $A_i$ is independent of $i$. Indeed, $A_i$ is built out of $h\in B^1(X)$ and $c\in B^2(X)$ (Corollary \ref{Fr2}), and so (since the number of possibilities for the coefficients of $A$ is countably infinite, whereas the base $\MM_{GM4}$ is uncountable) up to further shrinking the base $\MM_{GM4}$ we can assume $A_i=:A$ is constant. For the same reason, $\alpha_i$ is independent of $i$. Pushing forward under $\nu\times\nu$, the above equalities then give
    \[ \begin{split} \alpha  (\Delta_X-\Gamma_\iota) &=  {1\over d}\  (\nu\times\nu)_\ast (\sum_{i=1}^d I_i)\cdot A +  {1\over d}\ (\nu\times\nu)_\ast(\sum_{i=1}^d B_i)\\
                                                                                 &=: I\cdot A + B  \ \ \ \ \ \  \hbox{in}\ B^4(X\times X)\ ,\\
                                                                                 \end{split}   \] 
                        where the correspondences $I$, $A$ and $B$ are now generically defined {\em with respect to $\MM_{dEPW}$} and thus satisfy all the requirements of \eqref{onlyalg}.
      The constant $\alpha\in\QQ$ being non-zero (Theorem \ref{zh}), this proves equality \eqref{onlyalg}, and hence Proposition \ref{heart}. This closes the proof of the theorem.   
         \end{proof}

\begin{remark}\label{rem} As the reader will have noticed, the proof of Theorem \ref{main} makes a detour via algebraic equivalence. This is necessary to make the step from (cycles generically defined over) the moduli space $\MM_{GM4}$ to (cycles generically defined over) the moduli space $\MM_{dEPW}$. Indeed, for this step we rely on the Franchettina property (Corollary \ref{Fr2}), which we do not know how to prove (and which is probably false) on the level of rational equivalence.

As long as we work with cycles generically defined over $\MM_{GM4}$, this detour via algebraic equivalence is not necessary. That is, the proof of Theorem \ref{main} for double EPW sextics that are attained by the Iliev--Manivel construction is actually somewhat simpler, in that the passage to cycles modulo algebraic equivalence is not needed.
\end{remark}

\medskip    
 \vskip1cm
\begin{nonumberingt} Thanks to the referee for helpful comments. Thanks to Kai and Len for enjoying many Moomin episodes with me.
\end{nonumberingt}

\newpage

%\author[Robert Laterveer]
%{Robert Laterveer}

%\address{Institut de Recherche Math\'ematique Avanc\'ee,
%CNRS -- Universit\'e 
%de Strasbourg,\
%7 Rue Ren\'e Des\-car\-tes, 67084 Strasbourg CEDEX,
%FRANCE.}
%\email{robert.laterveer@math.unistra.fr}
%

\section{Erratum to ``The Beauville--Voisin conjecture for double EPW sextics''}

%\begin{abstract}  \end{abstract}

%\thanks{\textit{2020 Mathematics Subject Classification:}  14C15, 14C25, 14C30}
%\keywords{Algebraic cycles, Chow group, motive, hyperk\"ahler varieties, Beauville--Voisin conjecture, Beauville's ``splitting property'' conjecture, }
%\thanks{Supported by ANR grant ANR-20-CE40-0023.}

%\keywords{Algebraic cycles, Chow group, motive, Beauville's ``splitting property'' conjecture, multiplicative Chow--K\"unneth decomposition, Fano threefold, homological projective duality}
%\subjclass[2010]{Primary 14C15, 14C25, 14C30.}

%\maketitle

The last part of Remark 2.26 in \cite{BVEPW} is incorrect, and hence \cite[Equation (3)]{BVEPW} is not justified. This affects the proof of \cite[Proposition 2.27]{BVEPW}. To remedy this, we propose another definition of the cycle $I_\pi$; this definition replaces \cite[Definition 2.25]{BVEPW}.  The results in \cite{BVEPW} (apart from \cite[Remark 2.26]{BVEPW})
remain valid, provided one works everywhere with this new definition of the cycle $I_\pi$.

   \begin{definition}\label{def} Let $Y$ be a Gushel--Mukai fourfold that is general in the sense of \cite[Theorem 2.7]{BVEPW}, let $F=F(Y)$ be the variety of conics, and let $\pi\colon F\to X$ be the morphism to the associated
   double EPW sextic. We define the cycle $I_\pi\in A^2 (X\times X)$ by
     \[  I_\pi:= {1\over \delta}(\pi\times\pi)_\ast \bigl( I_F\cdot (h_F\times h_F)\bigr)\ \ \ \in\ A^2(X\times X)\ ,\]
     where $h_F\in A^1(F)$ is an ample divisor, and $\delta\in\NN$ is the intersection number of $h_F\times h_F$ with a general fiber of $\pi\times\pi$.
      \end{definition} 
      
%      \begin{remark} We note that a priori, the subvariety $I_\pi$ is {\em not \/} intrinsic to $X$. That is, if $Y_1, Y_2$ are two Gushel--Mukai fourfolds with the same associated double EPW sextic $X$
%      (such $Y_1$ and $Y_2$ are called ``period partners''), we do not know whether $I_{\pi_1}$ and $I_{\pi_2}$ coincide.
%      
%       We also remark that if $(c,c^\prime)\in I_F$ are intersecting conics that are both $\tau$-conics, then $I_F$ contains the whole fiber $(\pi\times\pi)^{-1}(\pi(c),\pi(c^\prime))$.
%      Thus, we have a relation
%      \begin{equation}\label{fiber}   I_F = (\pi\times\pi)^\ast (I_\pi) + R\ \ \ \hbox{in}\ A^2(F\times F)\ ,\end{equation}
%      where $R$ is supported on $\bigl((F_\rho(Y)\cup F_\sigma(Y))\times F\bigr) \cup \bigl(F\times (F_\rho(Y)\cup F_\sigma(Y))\bigr)$.      
%        \end{remark}           
% 

We establish a preliminary lemma, that replaces \cite[Equation (3)]{BVEPW}:

\begin{lemma}\label{lemma} Notation as above. Then
  \[  I_F - (\pi\times\pi)^\ast (I_\pi) = A + C \ \ \hbox{in}\ A^2(F\times F)\ ,\]
  where $A$ is decomposable and generically defined, i.e.
  \[      A\ \ \in\ \Bigl\langle   (p_i)^\ast  GDA^\ast_{\MM_{GM4}}(F)\Bigr\rangle\ ,\] 
  and $C$ is supported on $\bigl((F_\rho(Y)\cup F_\sigma(Y))\times F\bigr) \cup \bigl(F\times (F_\rho(Y)\cup F_\sigma(Y))\bigr)$.  
    \end{lemma}
  
  \begin{proof} Let us write $\omega$ for the difference
    \[ \omega:=  I_F - (\pi\times\pi)^\ast (I_\pi)   \ \ \hbox{in}\ A^2(F\times F)\ .\]  
  Let
    \[ \pi^\circ\colon F^\circ:=F_\tau(Y)\to X^\circ:=X\setminus (  x_\rho\cup x_\sigma) \]
    denote the restriction of $\pi$ to the dense open $F^\circ$, so that $\pi^\circ$ is a $\PP^1$-bundle \cite[Proposition 4.11]{IM}.
    
   The projective bundle formula for $\pi^\circ\times\pi^\circ$ tells us that
   \[ \begin{split} \omega\vert_{F^\circ\times F^\circ} = (\pi^\circ\times\pi^\circ)^\ast (a_2) + (\pi^\circ\times\pi^\circ)^\ast(a_1^1)\cdot (h_{F^\circ}\times F^\circ) +   (\pi^\circ\times\pi^\circ)^\ast(a_1^2)\cdot (F^\circ\times h_{F^\circ})&\\    + e (h_{F^\circ}\times h_{F^\circ})&\ ,\\
   \end{split}\]
   where $e\in\QQ$ and $a_2\in A^2(X^\circ\times X^\circ), a_1^i\in A^1(X^\circ\times X^\circ)$, and $h_{F^\circ}$ denotes the restriction of $h_F$ to $F^\circ$.
   
   Extending to $F\times F$, we thus find that
    \begin{equation}\label{this} \begin{split} \omega_{} = (\pi\times\pi)^\ast (\bar{a_2}) + (\pi\times\pi)^\ast(\bar{a_1^1})\cdot (h_F\times F) +   (\pi\times\pi)^\ast(\bar{a_1^2})\cdot (F\times h_F)&\\    + e (h_F\times h_F)+C&\ ,\\
   \end{split}\end{equation}
   where  $\bar{a_2}\in A^2(X\times X), \bar{a_1^i}\in A^1(X\times X)$, and $C$ is supported on  $\bigl((F_\rho(Y)\cup F_\sigma(Y))\times F\bigr) \cup \bigl(F\times (F_\rho(Y)\cup F_\sigma(Y))\bigr)$.

  We observe that the cycles $\bar{a_2}, \bar{a_1^i}$ are generically defined (with respect to $\MM_{GM4}$). The cycles $\bar{a_1^i}$, being divisors, are also decomposable. It follows that all but possibly the first term on the right-hand side of \eqref{this} are of the required form.
  
It remains to analyze the first term on the right-hand side of \eqref{this}. The projective bundle formula for $\pi^\circ\times\pi^\circ$ gives us that 
 \[ \begin{split} a_2 = &(\pi^\circ\times\pi^\circ)_\ast \bigl( \omega\cdot (h_{F^\circ}\times h_{F^\circ})\bigr)\\
 & + (c_1\times X^\circ)\cdot (\pi^\circ\times\pi^\circ)_\ast \bigl( \omega\cdot (h_{F^\circ}\times F^\circ)\bigr) \\
 & + (X^\circ\times c_1)\cdot (\pi^\circ\times\pi^\circ)_\ast \bigl( \omega\cdot (F^\circ\times h_{F^\circ})\bigr) \\
 &+ (c_1\times c_1)\cdot  (\pi^\circ\times\pi^\circ)_\ast(\omega)\ ,\\
 \end{split}\]
 where $c_1\in A^1(X^\circ)$ is the first Chern class of the vector bundle inducing $\pi^\circ$. 
 By definition of $I_\pi$, the first term is zero. 
 The last 3 terms are intersections of divisors on $X^\circ\times X^\circ$, hence decomposable. Since $c_1$ is generically defined, it follows that $a_2$ (and hence also $\bar{a_2}$) is generically defined and decomposable. That is, the right-hand side of \eqref{this} is of the required form $A+C$, proving the lemma.
  \end{proof}

    The following proposition replaces Proposition 2.27 in \cite{BVEPW}:
    
    \begin{proposition}\label{zh2} Let 
    $Y$ be a Gushel--Mukai fourfold that is general in the sense of \cite[Theorem 2.7]{BVEPW}, let $X$ be the associated double EPW sextic, and let $\iota\in \aut(X)$ denote the covering involution. There is a relation
      \[ (I_\pi)^2 = \alpha \Gamma_\iota + I_\pi\cdot A + B\ \ \ \hbox{in}\ A^4(X\times X)\ ,\]
      where $\alpha\in\QQ^\ast$, and $A$ and $B$ 
      are decomposable and generically defined.
     \end{proposition}
     
    \begin{proof} 
      Intersecting Zhang's quadratic relation on $F\times F$ \cite[Proposition 2.24]{BVEPW} with $h_F\times h_F$ and pushing forward to $X\times X$, we obtain
      \[   (\pi\times\pi)_\ast \Bigl( I_F^2\cdot (h_F\times h_F)   \Bigr) =     (\pi\times\pi)_\ast \Bigl( ( \alpha W + I_F\cdot A + B + C)\cdot (h_F\times h_F)\Bigr)\ \ \ \hbox{in}\ A^4(X\times X)\ .\]    
    Applying Lemma \ref{lemma} to the left-hand side, we obtain
     \[ \begin{split} (\pi\times\pi)_\ast \Bigl( I_F\cdot (h_F\times h_F) \cdot ( (\pi\times\pi)^\ast(I_\pi) + A^\prime+C^\prime) \Bigr) =&\\     (\pi\times\pi)_\ast \Bigl( ( \alpha W + I_F\cdot A + B +& C)\cdot (h_F\times h_F)\Bigr)
     \ \ \ \hbox{in}\ A^4(X\times X)\ ,\\
     \end{split}\]      
     where $A^\prime$ is decomposable and generically defined, and $C^\prime$ is supported on  $\bigl((F_\rho(Y)\cup F_\sigma(Y))\times F\bigr) \cup \bigl(F\times (F_\rho(Y)\cup F_\sigma(Y))\bigr)$.  
     
   This simplifies to
     \[  (\pi\times\pi)_\ast \Bigl( I_F\cdot (h_F\times h_F) \cdot  (\pi\times\pi)^\ast(I_\pi) \Bigr) =     (\pi\times\pi)_\ast \Bigl( ( \alpha W + I_F\cdot A^{\prime\prime} + B + C^{\prime\prime})\cdot (h_F\times h_F)\Bigr)\ \ \ \hbox{in}\ A^4(X\times X)\ ,\]     
     for some new $A^{\prime\prime}, C^{\prime\prime}$ that are still of the required form.
     
    Applying the projection formula and Definition \ref{def} to the left-hand side, we obtain
    \[ \delta (I_\pi)^2 =  (\pi\times\pi)_\ast \Bigl( ( \alpha W + I_F\cdot A^{\prime\prime} + B + C^{\prime\prime})\cdot (h_F\times h_F)\Bigr)\ \ \ \hbox{in}\ A^4(X\times X)\ ,\]
    where $\delta\in\NN^\ast$ is as in Definition \ref{def}.     
      
   From here on, the proof goes exactly as \cite[Proof of Proposition 2.27]{BVEPW}.   
                \end{proof} 

\medskip    
 \vskip1cm
\begin{nonumberingt} Thanks to Carlo Mazzanti for noticing the problem with \cite[Remark 2.26]{BVEPW}, and for asking me lots of questions about \cite{BVEPW}.

\end{nonumberingt}

\end{document}